\documentclass[12pt,a4paper,twoside]{article}
\usepackage[T2A]{fontenc}
\usepackage[english,russian]{babel}

\fontsize{12}{21} \selectfont

\newtheorem{theorem}{Theorem}[section]

\newtheorem{definition}{Definition}[section]

\begin{document}

\noindent

\makeatletter
\renewcommand{\@evenhead}{\hfil
{\bf Andrei I. Bodrenko.}}
\renewcommand{\@oddhead}{\hfil
\small{ \underline{The Christoffel problem and two analogs of the
Minkowski problem in Riemannian space.}} }

\noindent
\bigskip
\begin{center}
{\large \bf The Christoffel problem and two analogs of the
Minkowski problem in Riemannian space.}
\end{center}
\medskip

\begin{center}
{\bf Andrei~I.~Bodrenko} \footnote{\copyright  Andrei~I.~Bodrenko,
associate professor,
Department of Mathematics,\\
Volgograd State University,
University Prospekt 100, Volgograd, 400062, RUSSIA.\\
E.-mail: bodrenko@mail.ru \qquad \qquad http://www.bodrenko.com}
\end{center}

\begin{center}
{\bf Abstract}\\
\end{center}

{\small Author finds the solutions of the Christoffel problem for
open and closed surfaces in Riemannian space. The Christoffel
problem is reduced to the problem of construction the continuous
\\
G-deformations preserving the sum of principal radii of curvature
for every point of surface in Riemannian space. G-deformation
transfers every normal vector of surface in parallel along the
path of the translation for each point of surface. The following
analogs of the Minkowski problem for open and closed surfaces in
Riemannian space are being considered in this article: 1) the
problem of construction the surface with prescribed mean curvature
and condition of G-deformation; 2) the problem of construction the
deformations preserving the area of each arbitrary region of
surface and condition of G-deformation.}

\section*{Introduction}

\bigskip

The Christoffel problem (ChP) is well known fundamental problem of
 differential geometry. Author solves the ChP in Riemannian space
  as the problem of  finding the continuous $G-$deformations
  with prescribed the sum of principal radii of curvature.

In the article, there is being considered the problem of
construction the surface with prescribed mean curvature and
condition of G-deformation in Riemannian space, which is the
analog of the Minkowski problem.

 The second analog of the Minkowski problem is finding
the deformations preserving the area of each arbitrary region of
surface with condition of $G-$deformation.

Theorems 1 and 2 represents the properties of solutions of
considered problems for open and closed surfaces in Riemannian
space respectively.

\section*{\S 1.1. Basic definitions. Statement of the main results
for open surfaces in Riemannian space.}

Let $R^{3}$ be  the three-dimensional Riemannian space with metric
tensor $\tilde{a}_{\alpha\beta},$
  $F^{+}$ be the two-dimensional simply connected oriented
surface in $R^{3}$ with the boundary $\partial F.$

Let $F^{+}\in C^{m,\nu}, \nu \in (0;1) , m\ge 4.$ $\partial F\in
C^{m+1,\nu}.$ Let $F^{+}$ has all strictly positive principal
curvatures $k_{1}$ and $k_{2}.$ Let $F^{+}$ be oriented so that
mean curvature $H$ is strictly positive. Denote $K=k_{1}k_{2}.$

Let $F^{+}$ be given by immersion of the domain $D\subset E^{2}$
into $R^{3}$ by the equation: $y^{\sigma}=f^{\sigma+}(x), x\in D$,
$f:D \rightarrow R^{3}.$ Denote by $d\sigma(x)=\sqrt{g}
dx^{1}\wedge dx^{2}$ the area element of the surface $F^{+}$. We
identify the points of immersion of surface $F^{+}$ with the
corresponding coordinate sets in $R^{3}$. Without loss of
generality we assume that $D$ is unit disk. Let $x^{1}, x^{2}$ be
the Cartesian coordinates.

Symbol $_{,i}$ denotes covariant derivative in metric of surface
$F^{+}.$ Symbol $\partial_{i}$ denotes partial derivative by
variable $x^{i}.$ We will assume $\dot{f}\equiv \frac{d f}{d t}.$
We define $\Delta(f)\equiv f(t)-f(0).$

We consider continuous deformation of the surface $F^{+}$:
$\{F_{t}\}$ defined by the equations
$$y^{\sigma}_{t}=y^{\sigma} + z^{\sigma}(t), z^{\sigma}(0)\equiv 0,
 t\in[0;t_{0}], t_{0}>0. \eqno(1.1)$$

\begin{definition}
\label{definition 1}. Deformation $\{F_{t}\}$ is called the
continuous deformation preserving the sum of principal radii of
curvature ( or $Ch-$deformation ) if the following condition
holds: \\ $\Delta(\frac{1}{k_{1}}+\frac{1}{k_{2}})=0$ and
$z^{\sigma}(t)$ is continuous by $t,$
 where $k_{1}$ and $k_{2}$ are principal curvatures of $F^{+}$.

\end{definition}

\begin{definition}
\label{definition 2}. Deformation $\{F_{t}\}$ is called the
continuous deformation preserving the mean curvature ( or
$H-$deformation) if the following condition holds: $\Delta(H)=0$
and $z^{\sigma}(t)$ is continuous by $t.$

\end{definition}

\begin{definition}
\label{definition 3}.
Deformation $\{F_{t}\}$ is called the continuous $A-$deformation \\
if the following condition holds: $d\sigma_{t}-d\sigma=0$ and
$z^{\sigma}(t)$ is continuous by $t.$
\end{definition}

This means that $A-$deformation preserves the area of each
arbitrary region of surface.

The deformation $\{F_{t}\}$ generates the following set of paths in $R^{3}$
$$u^{\alpha_{0}}(\tau)=(y^{\alpha_{0}}+z^{\alpha_{0}}(\tau)), \eqno(1.2)$$
where $z^{\alpha_{0}}(0)\equiv 0, \tau \in [0;t], t\in[0;t_{0}], t_{0}>0.$

\begin{definition}
\label{definition 4}. The deformation $\{F_{t}\}$ is called the
$G-$deformation if every normal vector of surface transfers in
parallel along the path of the translation for each point of
surface.
\end{definition}

Indices denoted by Greek alphabet letters define tensor
coordinates in Riemannian space $R^{3}.$ We use the following
rule: a formula is valid for all admissible values of indices if
there are no instructions for which values of indices it is valid.
We use the Einstein rule. Let $g_{ij}$ and $b_{ij}$ be the
coefficients of the first and the second fundamental form
respectively.

Let, along the $\partial F$, be given vector field tangent to
$F^{+}.$ We denote it by the following formula:
$$v^{\alpha}=l^{i}y^{\alpha}_{,i}.  \eqno(1.3)$$

We consider the boundary-value condition:
$$\tilde{a}_{\alpha\beta}z^{\alpha}v^{\beta}=\tilde{\gamma}(s,t) ,
s\in \partial D. \eqno(1.4)$$

Let $v^{\alpha}$ and $\tilde{\gamma}$ be of class $C^{m-2,\nu}.$

We denote:
$$\tilde{\lambda}_{k}=\tilde{a}_{\alpha\beta}y^{\alpha}_{,k}v^{\beta},
k=1,2.  \eqno(1.5)$$
$$\lambda_{k}=\frac{ \tilde{\lambda}_{k}}
{(\tilde{\lambda}_{1})^{2}+(\tilde{\lambda}_{2})^{2}},
k=1,2.  \eqno(1.6)$$

$$\lambda(s)=\lambda_{1}(s)+i \lambda_{2}(s), s\in \partial D. \eqno(1.7)$$

Let $n$ be the index of the given boundary-value condition
$$n=\frac{1}{2\pi}\Delta_{\partial D} \arg \lambda(s).  \eqno(1.8)$$

\begin{theorem}
\label{theorem1}. Let $F^{+}\in C^{m,\nu}, \nu \in (0;1) , m\ge
4,$ $\partial F \in C^{m+1,\nu}.$ Let $\tilde{a}_{\alpha\beta} \in
C^{m,\nu},$ $\exists M_{0}=const>0$ such that
$\|\tilde{a}_{\alpha\beta}\|_{m,\nu}<M_{0},$ $\|\partial
\tilde{a}_{\alpha\beta}\|_{m,\nu}<M_{0},$ $\|\partial^{2}
\tilde{a}_{\alpha\beta}\|_{m,\nu}<M_{0}.$ Let $v^{\beta},
\tilde{\gamma} \in C^{m-2,\nu}(\partial D),$ $\tilde{\gamma}$ is
continuously differentiable by $t.$ Let, at the point
$(x^{1}_{(0)},x^{2}_{(0)})$ of the domain $D,$ the following
condition holds: $\forall t : z^{\sigma}(t)\equiv 0.$

Then the following statements hold:

1) if $n > 0$ then there exist $t_{0}>0$
and $\varepsilon(t_{0})>0$ such that
for any admissible $\tilde{\gamma}$ satisfying the condition:
$|\dot{\tilde{\gamma}}|_{m-2,\nu}\leq \varepsilon$
 for all $t\in [0, t_{0}):$

1a) there exists $(2n-1)-$parametric $ChG-$deformation of class
$C^{m-2,\nu}(\bar{D})$ continuous by $t.$

1b) there exists $(2n-1)-$parametric $HG-$deformation of class
$C^{m-2,\nu}(\bar{D})$ continuous by $t.$

1c) there exists $(2n-1)-$parametric $AG-$deformation of class
$C^{m-2,\nu}(\bar{D})$ continuous by $t.$

2) if $n < 0$ then there exist $t_{0}>0$
and $\varepsilon(t_{0})>0$ such that
for any admissible $\tilde{\gamma}$ satisfying the condition:
$|\dot{\tilde{\gamma}}|_{m-2,\nu}\leq \varepsilon(t_{0})$
 for all $t\in [0, t_{0}):$

2a) there exists at most one $ChG-$deformation of class
$C^{m-2,\nu}(\bar{D})$ continuous by $t.$

2b) there exists at most one $HG-$deformation of class
$C^{m-2,\nu}(\bar{D})$ continuous by $t.$

2c) there exists at most one $AG-$deformation of class
$C^{m-2,\nu}(\bar{D})$ continuous by $t.$

3) if $n = 0$ then there exist $t_{0}>0$
and $\varepsilon(t_{0})>0$ such that
for any admissible $\tilde{\gamma}$ satisfying the condition:
$|\dot{\tilde{\gamma}}|_{m-2,\nu}\leq \varepsilon$
 for all $t\in [0, t_{0}):$

3a) there exists one $ChG-$deformation of class
$C^{m-2,\nu}(\bar{D})$ continuous by $t.$

3b) there exists one $HG-$deformation of class
$C^{m-2,\nu}(\bar{D})$ continuous by $t.$

3c) there exists one $AG-$deformation of class
$C^{m-2,\nu}(\bar{D})$ continuous by $t.$

\end{theorem}

\section*{\S 1.2. Statement of the main results
for closed surfaces in Riemannian space.}

Let $F$ be the two-dimensional simply connected oriented closed
surface in $R^{3}.$

Let $F\in C^{m,\nu}, \nu \in (0;1) , m\ge 4.$  Let $F$ has all
strictly positive principal curvatures $k_{1}$ and $k_{2}$. Let F
be oriented so that mean curvature $H$ is strictly positive.

Let $F$ be glued from the two-dimensional simply connected
oriented surfaces $F^{+}$ and $F^{-}$ of class $C^{m,\nu}.$ Let
$F^{+}$ be attached to $F^{-}$ along the common boundary $\partial
F$ of class $C^{m+1,\nu}.$

Let $F^{+}$ and $F^{-}$ be given by immersions of the domain
$D\subset E^{2}$ into $R^{3}$ by the equation:
$y^{\sigma}=f^{\sigma\pm}(x), x\in D$, $f^{\pm}:D \rightarrow
R^{3}.$

\begin{theorem}
\label{theorem2}. Let $F\in C^{m,\nu}, \nu \in (0;1) , m\ge 4,$ be
closed surface. Let $F$ be glued from the two-dimensional simply
connected oriented surfaces $F^{+}$ and $F^{-}$ of class
$C^{m,\nu}.$ Let $F^{+}$ be attached to $F^{-}$ along the common
boundary $\partial F$ of class $C^{m+1,\nu}.$ Let
$\tilde{a}_{\alpha\beta} \in C^{m,\nu},$ $\exists M_{0}=const>0$
such that $\|\tilde{a}_{\alpha\beta}\|_{m,\nu}<M_{0},$ $\|\partial
\tilde{a}_{\alpha\beta}\|_{m,\nu}<M_{0},$ $\|\partial^{2}
\tilde{a}_{\alpha\beta}\|_{m,\nu}<M_{0}.$

1) Then there exists $t_{0}>0$ such that
 for all $t\in [0, t_{0}):$

1a) there exists three-parametric $ChG-$deformation of class
$C^{m-2,\nu}$ continuous by $t.$

1b) there exists three-parametric $HG-$deformation of class
$C^{m-2,\nu}$ continuous by $t.$

1c) there exists three-parametric $AG-$deformation of class
$C^{m-2,\nu}$ continuous by $t.$

2) If, at the point $T_{0}\in F^{+},$ the following additional
condition holds: $\forall t : z^{\sigma}(t)\equiv 0.$ Then there
exists $t_{0}>0$ such that
 for all $t\in [0, t_{0}):$

2a) there exists  only zero $ChG-$deformation of class
$C^{m-2,\nu}$ continuous by $t.$

2b) there exists  only zero $HG-$deformation of class
$C^{m-2,\nu}$ continuous by $t.$

2c) there exists  only zero $AG-$deformation of class
$C^{m-2,\nu}$ continuous by $t.$

\end{theorem}

We use all designations from [30, 32].

\section*{\S 2. Deduction the formulas of $ChG-$deformations,\\
$HG-$deformations and $AG-$deformations for surfaces in Riemannian
space.}

\section*{\S 2.1. The formulas of $G-$deformations, $\Delta(g)$ and $\Delta(k_{1}k_{2})$.}

We denote:
$$z^{\sigma}(t)=a^{j}(t)y^{\sigma},_{j} + c(t)n^{\sigma},  \eqno(2.1.1)$$
where $a^{j}(0)\equiv 0, c(0)\equiv 0,$ $n^{\sigma}$ is unit
normal vector of surface at the point $(y^{\sigma}).$ Therefore
the deformation of surface is defined by functions $a^{j}$ and
$c.$ We introduce conjugate isothermal coordinate system where
$b_{ii}=V,i=1,2, b_{12}=b_{21}=0.$

 The equations of $G-$deformation were obtained in [30] and [32]:

$$ \partial_{2}\dot{a}^{1}-\partial_{1}\dot{a}^{2}+p_{k}\dot{a}^{k}=
\dot{\Psi_{1}},\eqno(2.1.2) $$
 where $p_{k}$ and $\dot{\Psi}_{1}$ are defined in [30].
 Note that $p_{k}$ do not depend on $t.$

The function $\dot{c}$ is found on functions $\dot{a}^{i}$ from
formulas obtained in [30] and [32].

We have from [30]:
$$\Delta(g)=2g(
\partial_{1}a^{1}+\partial_{2}a^{2}+q_{k}a^{k}-\Psi_{2}), \eqno(2.1.3)$$
where $$q_{1}=\partial_{1}(\ln \sqrt{g}), q_{2}=\partial_{2}(\ln
\sqrt{g}), $$ Where $\Psi_{2}$ has explicit form and is defined in
[30]. Note that $q_{k}$ do not depend on $t.$

We obtain the following equation from [30]:\\
$$\Delta(K)=\frac{1}{b(t)}(
g\partial_{1}a^{1}+g\partial_{2}a^{2}+2gq_{k}a^{k}-2g\Psi_{2}-
\frac{g}{V}(M^{4}_{11}+M^{4}_{22})-\frac{g}{V^{2}}W_{2}^{(b)}),
\eqno(2.1.4) $$ Where $M^{4}_{11}, M^{4}_{22}, W_{2}^{(b)}$ have
explicit forms and are defined in [30].

Therefore we obtain:
$$\dot{\Delta}(K)=\frac{g}{b}(
\partial_{1}\dot{a}^{1}+\partial_{2}\dot{a}^{2}+q^{(b)}_{k}\dot{a}^{k}-
\dot{\Psi}_{2}^{(b)}),  \eqno(2.1.5)$$ where
$\dot{\Psi}_{2}^{(b)}=q^{(b)}_{0}\dot{c}- P_{0}(\dot{a}^{1},
\dot{a}^{2},\partial_{i}\dot{a}^{j}).$ $P_{0}$ has explicit form.
Notice that $q^{(b)}_{k}\in C^{m-3,\nu},$ $q^{(b)}_{0} \in
C^{m-3,\nu}$ and
 do not depend on $t.$

{\bf Lemma 2.1.1.} {\it Let the following conditions hold:

1) metric tensor in $R^{3}$ satisfies the conditions: $\exists
M_{0}=const>0$ such that
$\|\tilde{a}_{\alpha\beta}\|_{m,\nu}<M_{0},$ $\|\partial
\tilde{a}_{\alpha\beta}\|_{m,\nu}<M_{0},$ $\|\partial^{2}
\tilde{a}_{\alpha\beta}\|_{m,\nu}<M_{0}.$

2) $\exists t_{0}>0$ such that $a^{k}(t), \partial_{i}a^{k}(t),
\dot{a}^{k}(t), \partial_{i}\dot{a}^{k}(t)$ are continuous by $t,
\forall t\in[0,t_{0}],$ $ a^{k}(0)\equiv 0,
\partial_{i}a^{k}(0)\equiv 0.$

3) $\exists t_{0}>0$ such that $a^{i}(t)\in C^{m-2,\nu} ,
\partial_{k} a^{i}(t)\in C^{m-3,\nu},$ $\forall t\in[0,t_{0}].$

Then $\exists t_{*}>0$ such that for all $t\in[0,t_{*})$ $P_{0}\in
C^{m-3,\nu}$ and the following inequality holds:
$$\|P_{0}(\dot{a}^{1}_{(1)},\dot{a}^{2}_{(1)})-
P_{0}(\dot{a}^{1}_{(2)},\dot{a}^{2}_{(2)})\|_{m-2,\nu}\leq
K_{9}(t)(\|\dot{a}^{1}_{(1)}-\dot{a}^{1}_{(2)}\|_{m-1,\nu}+
\|\dot{a}^{2}_{(1)}-\dot{a}^{2}_{(2)}\|_{m-1,\nu}),$$ where for
any $\varepsilon>0$ there exists $t_{0}>0$ such that
 for all $t\in[0,t_{0})$ the following inequality holds:
$K_{9}(t)<\varepsilon.$ }

The proof follows from [30].

\section*{\S 2.2. Deduction the formulas of $\Delta(H)$ and $\dot{\Delta}(H)$.}

The formula of mean curvature is:\\
$$2H=g^{ij}b_{ij}=g^{11}b_{11}+g^{22}b_{22},
2H(t)=g^{ij}(t)b_{ij}(t). \eqno(2.2.1) $$

Then we have:\\
$$2\Delta(H)=g^{ij}(t)b_{ij}(t)-g^{ij}b_{ij}. \eqno(2.2.2) $$

We use the following formulas:\\
$$g^{11}(t)=\frac{g_{22}(t)}{g(t)},g^{22}(t)=\frac{g_{11}(t)}{g(t)},
g^{12}(t)=g^{21}(t)=-\frac{g_{12}(t)}{g(t)}. \eqno(2.2.3) $$

Then we have:\\
$$\Delta(H)=\frac{1}{2g(t)}(
g_{22}(t)b_{11}(t)+g_{11}(t)b_{22}(t)
-g_{12}(t)(b_{12}(t)+b_{21}(t))-2g(t)H). \eqno(2.2.4) $$

Consider the following formula:
$$\Delta(H)=\frac{1}{2g(t)}(g_{22}(t)b_{11}(t)+g_{11}(t)b_{22}(t)
-g_{12}(t)(b_{12}(t)+b_{21}(t))-2g(t)H). \eqno(2.2.5) $$

Use the formulas:\\
$$g_{ij}(t)=g_{ij}+\Delta(g_{ij}), b_{ij}(t)=b_{ij}+\Delta(b_{ij}),
g(t)=g+\Delta(g). \eqno(2.2.6) $$

Then we obtain the equation:\\
$$\Delta(H)=\frac{1}{2g(t)}(
g_{22}b_{11}(t)+g_{11}b_{22}(t)-g_{12}(b_{12}(t)+b_{21}(t))+$$
$$\Delta(g_{22})b_{11}(t)+\Delta(g_{11})b_{22}(t)-
\Delta(g_{12})(b_{12}(t)+b_{21}(t)) -2gH-2\Delta(g)H).
\eqno(2.2.7) $$

Therefore we get the equation:\\
$$\Delta(H)=\frac{1}{2g(t)}(
g_{22}b_{11}+g_{11}b_{22}-g_{12}(b_{12}+b_{21})+
g_{22}\Delta(b_{11})+g_{11}\Delta(b_{22})-$$
$$g_{12}(\Delta(b_{12})+\Delta(b_{21}))+
\Delta(g_{22})b_{11}+\Delta(g_{11})b_{22}-
\Delta(g_{12})(b_{12}+b_{21})+$$
$$\Delta(g_{22})\Delta(b_{11})+\Delta(g_{11})\Delta(b_{22})-
\Delta(g_{12})(\Delta(b_{12})+\Delta(b_{21}))-2gH-2\Delta(g)H).
\eqno(2.2.8) $$

Simplifying we obtain the equation:\\
$$\Delta(H)=\frac{1}{2g(t)}(
g_{22}\Delta(b_{11})+g_{11}\Delta(b_{22})-
g_{12}(\Delta(b_{12})+\Delta(b_{21}))+
\Delta(g_{22})b_{11}+\Delta(g_{11})b_{22}-$$
$$\Delta(g_{22})\Delta(b_{11})+\Delta(g_{11})\Delta(b_{22})-
\Delta(g_{12})(\Delta(b_{12})+\Delta(b_{21}))-2\Delta(g)H).
\eqno(2.2.9) $$

Then we have:\\
$$\Delta(H)=\frac{1}{2g(t)}(
g_{22}\Delta(b_{11})+g_{11}\Delta(b_{22})-
g_{12}(\Delta(b_{12})+\Delta(b_{21}))+
V\Delta(g_{22})+V\Delta(g_{11})-$$
$$\Delta(g_{22})\Delta(b_{11})+\Delta(g_{11})\Delta(b_{22})-
\Delta(g_{12})(\Delta(b_{12})+\Delta(b_{21}))-2\Delta(g)H).
\eqno(2.2.10) $$

We use the formula:
$$\Delta(b_{ij})=\partial_{i}(a^{k})b_{jk}+M^{4}_{ij}. \eqno(2.2.11) $$
We can write the following:
$$\Delta(g_{ii})=\partial_{i}(a^{i})g_{ii}+M^{5}_{ii}. \eqno(2.2.12) $$

Therefore we have
$$\Delta(H)=\frac{1}{2g(t)}(
Vg_{22}\partial_{1}(a^{1})+Vg_{11}\partial_{2}(a^{2})+
Vg_{22}\partial_{2}(a^{2})+Vg_{11}\partial_{1}(a^{1})+$$
$$g_{22}M^{4}_{11}+g_{11}M^{4}_{22}+
VM^{5}_{22}+VM^{5}_{11} -g_{12}(\Delta(b_{12})+\Delta(b_{21}))$$
$$-\Delta(g_{22})\Delta(b_{11})+\Delta(g_{11})\Delta(b_{22})-
\Delta(g_{12})(\Delta(b_{12})+\Delta(b_{21}))-2\Delta(g)H).
\eqno(2.2.13) $$

Denote
$$\Psi_{4}=g_{22}M^{4}_{11}+g_{11}M^{4}_{22}+
VM^{5}_{22}+VM^{5}_{11} -g_{12}(\Delta(b_{12})+\Delta(b_{21}))$$
$$-\Delta(g_{22})\Delta(b_{11})+\Delta(g_{11})\Delta(b_{22})-
\Delta(g_{12})(\Delta(b_{12})+\Delta(b_{21})). \eqno(2.2.14) $$

Then we obtain the following equation
$$\Delta(H)=\frac{1}{2g(t)}(
Vg_{22}\partial_{1}(a^{1})+Vg_{11}\partial_{2}(a^{2})+
Vg_{22}\partial_{2}(a^{2})+Vg_{11}\partial_{1}(a^{1})+\Psi_{4}-2\Delta(g)H).
\eqno(2.2.15) $$

The equation takes the form
$$\Delta(H)=\frac{1}{2g(t)}(
V(g_{11}+g_{22})(\partial_{1}(a^{1})+\partial_{2}(a^{2}))+\Psi_{4}-
2\Delta(g)H). \eqno(2.2.16) $$

We use  formulas (2.1.20) and (2.1.21). Then we have
$$\Delta(H)=\frac{1}{2g(t)}(
V(g_{11}+g_{22})(\partial_{1}a^{1}+\partial_{2}a^{2})+\Psi_{4}-$$
$$4gH(\partial_{1}a^{1}+\partial_{2}a^{2}+q_{k}a^{k}-\Psi_{2})). \eqno(2.2.17) $$

Using the formula
$$2Hg=V(g_{11}+g_{22})$$
we get the following equation
$$\Delta(H)=\frac{1}{2g(t)}(
(-2gH(\partial_{1}a^{1}+\partial_{2}a^{2})-4Hgq_{k}a^{k}+
4Hg\Psi_{2}+\Psi_{4}). \eqno(2.2.18) $$

Therefore
$$\Delta(H)=\frac{Hg}{g(t)}(
-\partial_{1}a^{1}-\partial_{2}a^{2}-2q_{k}a^{k}+
2\Psi_{2}+\frac{\Psi_{4}}{2Hg}). \eqno(2.2.19) $$

Differentiating by $t$ we have
$$\dot{\Delta}(H)=\frac{Hg}{g(t)}(
-\partial_{1}\dot{a}^{1}-\partial_{2}\dot{a}^{2}-2q_{k}\dot{a}^{k}+
2\dot{\Psi_{2}}+\frac{\dot{\Psi}_{4}}{2Hg})-$$
$$\frac{Hg\dot{g}(t)}{(g(t))^2}(
-\partial_{1}a^{1}-\partial_{2}a^{2}-2q_{k}a^{k}+
2\Psi_{2}+\frac{\Psi_{4}}{2Hg}). \eqno(2.2.20)$$

Therefore we obtain:
$$\dot{\Delta}(H)=H(
-\partial_{1}\dot{a}^{1}-\partial_{2}\dot{a}^{2}-q^{(h)}_{k}\dot{a}^{k}+
\dot{\Psi}_{2}^{(h)}),  \eqno(2.2.21)$$ where
$\dot{\Psi}_{2}^{(h)}=q^{(h)}_{0}\dot{c}- P_{0}^{(h)}(\dot{a}^{1},
\dot{a}^{2},\partial_{i}\dot{a}^{j}).$ Notice that $q^{(h)}_{k}\in
C^{m-3,\nu},$ $q^{(h)}_{0} \in C^{m-3,\nu}$ and
 do not depend on $t.$

{\bf Lemma 2.2.1.} {\it Let the following conditions hold:

1) metric tensor in $R^{3}$ satisfies the conditions: $\exists
M_{0}=const>0$ such that
$\|\tilde{a}_{\alpha\beta}\|_{m,\nu}<M_{0},$ $\|\partial
\tilde{a}_{\alpha\beta}\|_{m,\nu}<M_{0},$ $\|\partial^{2}
\tilde{a}_{\alpha\beta}\|_{m,\nu}<M_{0}.$

2) $\exists t_{0}>0$ such that $a^{k}(t), \partial_{i}a^{k}(t),
\dot{a}^{k}(t), \partial_{i}\dot{a}^{k}(t)$ are continuous by $t,
\forall t\in[0,t_{0}],$ $ a^{k}(0)\equiv 0,
\partial_{i}a^{k}(0)\equiv 0.$

3) $\exists t_{0}>0$ such that $a^{i}(t)\in C^{m-2,\nu} ,
\partial_{k} a^{i}(t)\in C^{m-3,\nu},$ $\forall t\in[0,t_{0}].$

Then $\exists t_{*}>0$ such that for all $t\in[0,t_{*})$
$P_{0}^{(h)}\in C^{m-3,\nu}$ and the following inequality holds:
$$\|P_{0}^{(h)}(\dot{a}^{1}_{(1)},\dot{a}^{2}_{(1)})-
P_{0}^{(h)}(\dot{a}^{1}_{(2)},\dot{a}^{2}_{(2)})\|_{m-2,\nu}\leq
K_{10}(t)(\|\dot{a}^{1}_{(1)}-\dot{a}^{1}_{(2)}\|_{m-1,\nu}+
\|\dot{a}^{2}_{(1)}-\dot{a}^{2}_{(2)}\|_{m-1,\nu}),$$ where for
any $\varepsilon>0$ there exists $t_{0}>0$ such that
 for all $t\in[0,t_{0})$ the following inequality holds:
$K_{10}(t)<\varepsilon.$ }

The proof follows from construction of function $P_{0}^{(h)}$ and
lemmas of \S 7 and \S 8.

Notice the following formula:
$$\Delta(H)=\frac{1}{2g(t)}(g_{22}\Delta(b_{11})+g_{11}\Delta(b_{22})-$$
$$g_{12}(\Delta(b_{12})+\Delta(b_{21}))+
V\Delta(g_{22})+V\Delta(g_{11})-$$
$$\Delta(g_{22})\Delta(b_{11})+\Delta(g_{11})\Delta(b_{22})-
\Delta(g_{12})(\Delta(b_{12})+\Delta(b_{21}))-2\Delta(g)H).
\eqno(2.2.22) $$

Therefore we get the following formula
$$\Delta(H)=\frac{1}{2g(t)}
(-2gH(\partial_{1}a^{1}+\partial_{2}a^{2})-4Hgq_{k}a^{k}+
4Hg\Psi_{2}+\Psi_{4}). \eqno(2.2.23) $$

Therefore we obtain:\\

$$\dot{\Delta}(H)=-\frac{\dot{g}(t)}{2(g(t))^{2}}
(g_{22}\Delta(b_{11})+g_{11}\Delta(b_{22})-$$
$$g_{12}(\Delta(b_{12})+\Delta(b_{21}))+
V\Delta(g_{22})+V\Delta(g_{11})-$$
$$\Delta(g_{22})\Delta(b_{11})+\Delta(g_{11})\Delta(b_{22})-
\Delta(g_{12})(\Delta(b_{12})+\Delta(b_{21}))-2\Delta(g)H)+$$
$$\frac{1}{2g(t)}(g_{22}\dot{\Delta}(b_{11})+g_{11}\dot{\Delta}(b_{22})-
g_{12}(\dot{\Delta}(b_{12})+\dot{\Delta}(b_{21}))+
V\dot{\Delta}(g_{22})+V\dot{\Delta}(g_{11})-$$
$$\dot{\Delta}(g_{22})\Delta(b_{11})+\dot{\Delta}(g_{11})\Delta(b_{22})-
\dot{\Delta}(g_{12})(\Delta(b_{12})+\Delta(b_{21}))+$$
$$\Delta(g_{22})\dot{\Delta}(b_{11})+\Delta(g_{11})\dot{\Delta}(b_{22})-
\Delta(g_{12})(\dot{\Delta}(b_{12})+\dot{\Delta}(b_{21}))-
2\dot{\Delta}(g)H). \eqno(2.2.24) $$

Hence we can write the following
$$\dot{\Delta}(H)=\frac{1}{2g(t)}
(-2gH(\partial_{1}\dot{a}^{1}+\partial_{2}\dot{a}^{2})-4Hgq_{k}\dot{a}^{k}+
4Hg\dot{\Psi}_{2}+\dot{\Psi}_{4})$$
$$
-\frac{\dot{g}(t)}{2(g(t))^{2}}
(-2gH(\partial_{1}a^{1}+\partial_{2}a^{2})-4Hgq_{k}a^{k}+
4Hg\Psi_{2}+\Psi_{4}). \eqno(2.2.25) $$

We obtain
$$\Delta(H)=\frac{1}{2g(t)}
(-2gH(\partial_{1}a^{1}+\partial_{2}a^{2}+
q^{(h)}_{k}a^{k})+\Psi_{2}^{(h)}). \eqno(2.2.26) $$

Therefore we get
$$\dot{\Delta}(H)=\frac{1}{2g(t)}
(-2gH(\partial_{1}\dot{a}^{1}+\partial_{2}\dot{a}^{2}+
q^{(h)}_{k}\dot{a}^{k})+\dot{\Psi}_{2}^{(h)})$$
$$-\frac{\dot{g}(t)}{2(g(t))^{2}}
(-2gH(\partial_{1}a^{1}+\partial_{2}a^{2}+
q^{(h)}_{k}a^{k})+\Psi_{2}^{(h)}). \eqno(2.2.27) $$

\section*{\S 2.3. Deduction the formulas of deformations preserving
the sum of principal radii of curvature.}

We have the formula
$$\Delta(\frac{H}{K})=\frac{H(t)}{K(t)}-\frac{H}{K}. \eqno(2.3.1) $$

Therefore we obtain the equation of $Ch-$deformation
 preserving the sum of principal radii of curvature.
$$\Delta(H)=\frac{H}{K}\Delta(K). \eqno(2.3.2) $$

Using formulas from \S 2.1. and \S 2.2. we have
$$ \partial_{1}\dot{a}^{1}+\partial_{2}\dot{a}^{2}+q^{(c)}_{k}\dot{a}^{k}=
\dot{\Psi}^{(c)}_{2}, \eqno(2.3.3) $$ where
$\dot{\Psi}^{(c)}_{2}=q^{(c)}_{0}\dot{c}-P_{0}^{(c)}.$
 Note that $q^{(c)}_{k}$ do not depend on $t,$
$P_{0}^{(c)}(\dot{a}^{1}, \dot{a}^{2},\partial_{i}\dot{a}^{j}).$
Notice that $q^{(c)}_{k}\in C^{m-3,\nu},$ $q^{(c)}_{0} \in
C^{m-3,\nu}$ and
 do not depend on $t.$

{\bf Lemma 2.3.1.} {\it Let the following conditions hold:

1) metric tensor in $R^{3}$ satisfies the conditions: $\exists
M_{0}=const>0$ such that
$\|\tilde{a}_{\alpha\beta}\|_{m,\nu}<M_{0},$ $\|\partial
\tilde{a}_{\alpha\beta}\|_{m,\nu}<M_{0},$ $\|\partial^{2}
\tilde{a}_{\alpha\beta}\|_{m,\nu}<M_{0}.$

2) $\exists t_{0}>0$ such that $a^{k}(t), \partial_{i}a^{k}(t),
\dot{a}^{k}(t), \partial_{i}\dot{a}^{k}(t)$ are continuous by $t,
\forall t\in[0,t_{0}],$ $ a^{k}(0)\equiv 0,
\partial_{i}a^{k}(0)\equiv 0.$

3) $\exists t_{0}>0$ such that $a^{i}(t)\in C^{m-2,\nu} ,
\partial_{k} a^{i}(t)\in C^{m-3,\nu},$ $\forall t\in[0,t_{0}].$

Then $\exists t_{*}>0$ such that for all $t\in[0,t_{*})$
$P_{0}^{(c)}\in C^{m-3,\nu}$ and the following inequality holds:
$$\|P_{0}^{(c)}(\dot{a}^{1}_{(1)},\dot{a}^{2}_{(1)})-
P_{0}^{(c)}(\dot{a}^{1}_{(2)},\dot{a}^{2}_{(2)})\|_{m-2,\nu}\leq
K_{15}(t)(\|\dot{a}^{1}_{(1)}-\dot{a}^{1}_{(2)}\|_{m-1,\nu}+
\|\dot{a}^{2}_{(1)}-\dot{a}^{2}_{(2)}\|_{m-1,\nu}),$$ where for
any $\varepsilon>0$ there exists $t_{0}>0$ such that
 for all $t\in[0,t_{0})$ the following inequality holds:
$K_{15}(t)<\varepsilon.$ }

The proof follows from construction of function $P_{0}^{(c)}$ and
lemmas of \S 7 and \S 8 of [30].

The equation (2.3.3) determines deformations of surface
preserving the sum of \\ principal radii of curvature with
condition of $G-$deformation.

\section*{\S 3. Proof of theorems 1 and 2.}

We have the following equation systems of elliptic type

a) for $ChG-$deformations:
$$ \partial_{2}\dot{a}^{1}-\partial_{1}\dot{a}^{2}+p_{k}\dot{a}^{k}=
\dot{\Psi_{1}}, $$
$$ \partial_{1}\dot{a}^{1}+\partial_{2}\dot{a}^{2}+q^{(c)}_{k}\dot{a}^{k}=
\dot{\Psi}^{(c)}_{2}, \eqno(3.1a) $$ where we use (2.1.2) and
(2.3.3). $\dot{\Psi}^{(c)}_{2}=q^{(c)}_{0}\dot{c}-P_{0}^{(c)}.$
 Note that $q^{(c)}_{k}$ do not depend on $t.$

b) for $HG-$deformations:
$$ \partial_{2}\dot{a}^{1}-\partial_{1}\dot{a}^{2}+p_{k}\dot{a}^{k}=
\dot{\Psi_{1}}, $$
$$ \partial_{1}\dot{a}^{1}+\partial_{2}\dot{a}^{2}+q^{(h)}_{k}\dot{a}^{k}=
\dot{\Psi}^{(h)}_{2}, \eqno(3.1b) $$ where we use (2.1.2) and
(2.2.21). $\dot{\Psi}^{(h)}_{2}=q^{(h)}_{0}\dot{c}-P_{0}^{(h)}.$
 Note that $q^{(h)}_{k}$ do not depend on $t.$

c) for $AG-$deformations:
$$ \partial_{2}\dot{a}^{1}-\partial_{1}\dot{a}^{2}+p_{k}\dot{a}^{k}=
\dot{\Psi_{1}}, $$
$$ \partial_{1}\dot{a}^{1}+\partial_{2}\dot{a}^{2}+q_{k}\dot{a}^{k}=
\dot{\Psi}_{2}, \eqno(3.1c) $$ where we use (2.1.2) and (2.1.3).
$\dot{\Psi}_{2}$ is defined in [30].
 Note that $q_{k}$ do not depend on $t.$

For theorem 1, we reduce (3.1a), (3.1b) and (3.1c) with
boundary-value condition (1.4) by the methods form [30] to the
following form of desired boundary-value problem:

$$\partial_{\bar{z}}\dot{w}+A\dot{w}+B\bar{\dot{w}}+E(\dot{w})=\dot{\Psi},
\qquad Re\{\overline{\lambda}\dot{w}\}=\dot{\varphi} \quad on
\quad \partial D, \eqno(3.14) $$ where $\dot{\Psi},$
$\dot{\varphi},$ $E$ have explicit form and are defined in a
similar way as it was made in [30],
$\lambda=\lambda_{1}+i\lambda_{2},$ $|\lambda|\equiv 1,$ $\lambda,
\dot{\varphi}\in C^{m-2,\nu}(\partial D).$

We use estimations of norms for obtained functions from \S2, \S3,
\S7 of [30]. Formulas of functions $\dot{W}_{1},\dot{W}_{2},
\dot{\Psi}_{2},$ the estimations for norms of  these functions are
presented in \S8 of [30]. From article [30], we also use the
following lemmas: 2.1, 2.2, 2.3, 3.1, 3.2, 5.1, 6.2.1, 7.1, 7.2,
7.3, 8.3.1, 8.3.2, and theorems: 1, 9.1, 9.2.

Therefore proof of theorem 1 follows from using similar reasonings
for proving theorem 1 from article [30].

We obtain the proof of theorem 2 by using methods for proving
theorem 1 from article [32].

\section*{\bf References.}

\begin{enumerate}

 \item A.I. Bodrenko.
On continuous almost ARG-deformations of hypersurfaces in
Euclidean space [in Russian]. Dep. in VINITI 27.10.92., N3084-T92,
UDK 513.81, 14 pp.

 \item A.I. Bodrenko. Some properties continuous ARG-deformations
 [in Russian]. Theses of international science conference
"Lobachevskii and modern geometry", Kazan, Kazan university
publishing house, 1992 ., pp.15-16.

 \item A.I. Bodrenko. On continuous ARG-deformations [in Russian].
Theses of reports on republican science and methodical conference,
dedicated to the 200-th anniversary of N.I.Lobachevskii, Odessa,
Odessa university publishing house, 1992 ., Part 1, pp.56-57.

 \item A.I. Bodrenko. On extension of infinitesimal almost
 ARG-deformations closed \\
 hypersurfaces into analytic
 deformations in Euclidean spaces [in Russian].
 Dep. in VINITI 15.03.93., N2419-T93 UDK 513.81, 30 pp.

 \item A.I. Bodrenko. On extension of infinitesimal almost
 ARG-deformations of hypersurface with boundary into analytic deformations
[in Russian]. Collection works of young scholars of VolSU,
Volgograd, Volgograd State University publishing house, 1993,\\
pp.79-80.

 \item A.I. Bodrenko. Some properties of continuous almost
AR-deformations of hypersurfaces with prescribed change of
Grassmannian image  [in Russian].
 Collection of science works of young scholars, Taganrog,
 Taganrog State Pedagogical Institute publishing house , 1994, pp. 113-120.

 \item A.I. Bodrenko. On continuous almost AR-deformations
with prescribed change of Grassmannian image [in Russian].
All-Russian school-colloquium on stochastic \\ methods of geometry
and analysis. Abrau-Durso. Publisher Moscow: "TVP". Theses of
reports, 1994, pp. 15-16.

 \item A.I. Bodrenko. Extension of infinitesimal almost ARG-deformations
of \\ hypersurfaces into analytic deformations [in Russian].
All-Russian school-colloquium on stochastic methods. Yoshkar-Ola.
Publisher Moscow: "TVP". Theses of reports, 1995, pp. 24-25.

 \item A.I. Bodrenko. Areal-recurrent deformations of hypersurfaces
 preserving Grassmannian image [in Russian].
Dissertation of candidate of physical-mathematical sciences. \\
Novosibirsk,1995, pp. 85.

 \item A.I. Bodrenko. Areal-recurrent deformations of hypersurfaces
preserving Grassmannian image [in Russian]. Author's summary of
dissertation of candidate of physical- \\ mathematical sciences.
Novosibirsk, 1995, pp. 1-14.

 \item A.I. Bodrenko. Some properties of ARG-deformations [in Russian].
Izvestiay Vuzov. Ser. Mathematics, 1996, N2, pp.16-19.

 \item A.I. Bodrenko. Continuous almost ARG-deformations of surfaces
 with boundary [in Russian].
Modern geometry and theory of physical fields. \\ International
geometry seminar of N.I.Lobachevskii Theses of reports, Kazan,
\\
Publisher Kazan university, 1997, pp.20-21.

 \item A.I. Bodrenko. Continuous almost AR-deformations of surfaces
with prescribed change of Grassmannian image [in Russian]. Red.
"Sib. mat. zhurnal.", Sib. otd. RAN , Novosibirsk, Dep. in VINITI
13.04.98., N1075-T98 UDK 513.81, 13 pp.

 \item A.I. Bodrenko. Almost ARG-deformations of the second order of
 surfaces in \\ Riemannian space [in Russian].
 Surveys in Applied and Industrial Mathematics. \\ 1998,
 Vol. 5, Issue 2, p.202. Publisher Moscow: "TVP".

 \item A.I. Bodrenko. Almost $AR$-deformations of a surfaces with prescribed
change of \\ Grassmannian image with exterior connections [in
 Russian]. Red. zhurn. "Izvestya vuzov. Mathematics.", Kazan, Dep.
in VINITI  03.08.98, N2471 - B 98. P. 1-9.

 \item A.I. Bodrenko. Properties of generalized G-deformations with
areal condition of normal type in Riemannian space [in Russian].
 Surveys in Applied and Industrial Mathematics.
Vol. 7. Issue 2. (VII All-Russian school-colloquium on stochastic
methods. Theses of reports.) P. 478. Moscow: TVP, 2000.

 \item I.N. Vekua. Generalized Analytic Functions. Pergamon.
 New York. 1962.

 \item I.N. Vekua. Generalized Analytic Functions [in Russian].
 Moscow. Nauka. 1988.

 \item I.N. Vekua. Some questions of the theory of differential equations
 and  applications in mechanics [in Russian].
 Moscow:"Nauka". 1991 . pp. 256.

 \item A.V. Zabeglov. On decidability of one nonlinear
  boundary-value problem for AG-deformations
  of surfaces with boundary [in Russian].
  Collection of science works. \\ Transformations of surfaces,
  Riemannian spaces determined by given recurrent \\ relations. Part 1.
Taganrog. Taganrog State Pedagogical Institute publishing house.
1999. pp. 27-37.

 \item V.T. Fomenko. On solution of the generalized Minkowski problem
 for surface with boundary [in Russian].
Collection of science works. Transformations of surfaces,\\
  Riemannian spaces determined by given recurrent relations. Part 1.
Taganrog. \\ Taganrog State Pedagogical Institute publishing
house. 1999. pp. 56-65.

 \item V.T. Fomenko. On uniqueness of solution of the generalized
 Christoffel problem for surfaces with boundary [in Russian].
Collection of science works. Transformations of surfaces,
  Riemannian spaces determined by given recurrent relations. Part 1.
Taganrog. Taganrog State Pedagogical Institute publishing house.
1999. pp. 66-72.

 \item V.T. Fomenko. On rigidity of surfaces with boundary in
 Riemannian space [in Russian].
 Doklady Akad. Nauk SSSR. 1969 . Vol. 187, N 2, pp. 280-283.

 \item V.T. Fomenko. ARG-deformations of hypersurfaces in Riemannian
space [in Russian].\\  //Dep. in VINITI 16.11.90 N5805-B90

  \item S.B. Klimentov. On one method of construction the solutions
  of boundary-value \\ problems in the bending theory of surfaces of positive
  curvature [in Russian]. \\ Ukrainian geometry sbornik. pp. 56-82.

 \item M.A. Krasnoselskii. Topological methods in the theory of
 nonlinear problems

   [in Russian]. Moscow, 1965.

 \item L.P. Eisenhart. Riemannian geometry [in Russian]. \\ Izd. in. lit.,
 Moscow 1948.
 (Eisenhart Luther Pfahler. Riemannian geometry. 1926.)

 \item J.A. Schouten, D.J. Struik. Introduction into new methods
 of differential geometry [in Russian]. Volume 2. Moscow. GIIL. 1948 .
 (von J.A. Schouten und D.J. Struik.
 Einf$\ddot{u}$hrung in die neueren methoden der \\
 differentialgeometrie.
 Zweite vollst$\ddot{a}$ndig umgearbeitete
 Auflage. Zweiter band. 1938. )

 \item I. Kh. Sabitov. //VINITI. Results of science and technics.
 Modern problems of \\mathematics [in Russian].
 Fundamental directions. Vol.48, pp.196-271.

 \item Andrei I. Bodrenko.The solution of the Minkowski problem for
open surfaces in \\ Riemannian space. Preprint 2007.
arXiv:0708.3929v1 [math.DG]

\item F.D. Gakhov. Boundary-value problems [in Russian]. GIFML.
Moscow. 1963.

\item Andrei I. Bodrenko.The solution of the Minkowski problem for
closed surfaces in \\ Riemannian space. Preprint 2007.
 arXiv:math.DG

\end{enumerate}

\end{document}